\documentclass[x11names,final,hidelinks]{article}

\usepackage{xcolor}

\usepackage{amsmath}
\usepackage{amsfonts}
\usepackage{graphicx}
\usepackage{hyperref}
\usepackage{cleveref}
\usepackage{geometry}
\usepackage{authblk}
\usepackage{enumerate}

\crefformat{equation}{(#2#1#3)}
\crefrangeformat{equation}{(#3#1#4) to~(#5#2#6)}
\crefmultiformat{equation}{(#2#1#3)}{ and~(#2#1#3)}{, (#2#1#3)}{ and~(#2#1#3)}

\newcommand{\extfig}[1]{\includegraphics{img/#1.pdf}}

\title{The geometry of rest-spike bistability\thanks{GIC has received funding from EPSRC (RG80792) and Qualcomm Inc. The research leading to these results has received funding from the European Research Council under the Advanced ERC Grant Agreement Switchlet n.670645.}}

\author{Giuseppe Ilario Cirillo\thanks{corresponding author: gic27@cam.ac.uk} }
\author{Rodolphe Sepulchre}
\affil{Department of Engineering, University of Cambridge, Cambridge, U.K.}
\date{}

\ifpdf
\hypersetup{
  pdftitle={The geometry of rest-spike bistability},
  pdfauthor={G. I. Cirillo, R. Sepulchre}
}
\fi

\begin{document}

\maketitle

\begin{abstract}
Morris-Lecar model is arguably the simplest dynamical model that retains both the slow-fast geometry of excitable phase portraits and the physiological interpretation of a conductance-based model. We augment this model with one slow inward current to capture the additional property of bistability between a resting state and a spiking limit cycle for a range of input current. The resulting dynamical system is a core structure for many dynamical phenomena such as slow spiking and bursting. We show how the proposed model combines physiological interpretation and mathematical tractability  and we discuss the benefits of the proposed approach with respect to alternative models in the literature.
\end{abstract}

\section{Introduction}

Conductance-based models are by now well established as a fundamental modeling framework to connect the physiology and the dynamics of excitable cells. Ever since the seminal work of Hodgkin and Huxley~\cite{Hodgkin_1952}, there has been a continuing effort in the literature to develop models that combine mathematical tractability and physiological interpretation. An interesting example is the two-dimensional model published by Morris and Lecar  in 1981 \cite{Morris_1981}. Like Hodgkin-Huxley model, it captures the essential physiology of excitability: a spike results from the fast activation of an inward current followed by the slow activation of an outward current. The former provides positive feedback in the fast time-scale whereas the latter provides negative feedback in the slow time-scale. Because it is only two-dimensional, the model is also amenable to phase-portrait analysis without any reduction. Its geometry is similar to the one of FitzHugh-Nagumo model~\cite{FitzHugh_1961}, the first mathematical model proposed to understand the core dynamics of the Hodgkin-Huxley model. In that sense, Morris-Lecar model combines the physiological interpretation of Hodgkin-Huxley model and the mathematical tractability of Fitzugh-Nagumo circuit.

In the present paper, we aim at capturing in a similar way the essence of rest-spike bistability, that is, the coexistence of a stable spiking attractor and a stable fixed point in a slow-fast model. The importance of this phenomenon is well acknowledged in the neurodynamics literature due to its role as a building block of neuronal patterns such as bursting~\cite{Ermentrout_2010,Izhikevich_2006}. We obtain rest-spike bistability by adding one extra current in Morris-Lecar model: an inward current with slow activation. The resulting model combines the three following features:
\begin{enumerate}[i)]
  \item For a range of input currents, the model is rest-spike bistable, that is, a stable equilibrium coexists with a stable limit cycle. The geometry of the two attractors is robust to the time-scale separation in the sense that it persist in the limit of infinite time-scale separation.
  \item The model has the direct physiological interpretation of dynamics and attractors being shaped by three distinct currents (fast positive feedback (e.g. sodium activation), slow negative feedback (e.g. potassium activation), and slow positive feedback (e.g. calcium activation)).
  \item The model is amenable to a mathematical analysis by geometric singular perturbation theory.
\end{enumerate}
We are not aware of other single-cell models in the literature combining those three features. Mathematical models of rest-spike bistability often lack the first feature above. For instance, a homoclinic bifurcation in the Morris-Lecar model only exists for a specific time-scale separation (see e.g. Table 3.1 in section 3.2 of~\cite{Ermentrout_2010}). Limitations of such models with respect to the geometry of the attractors and the robustness of bursting are discussed in~\cite{Franci_2012,Franci_2018}. We are only aware of two published models in which rest-spike bistability persists in the singular limit of infinite time-scale separation. The first one is the model proposed by Hindmarsh and Rose in 1982~\cite{Hindmarsh_1982} as a mathematical model aimed at capturing low-frequency spiking. The second one is the transcritical model proposed in 2011 as a two-dimensional reduction of a physiological model combining the currents of the Hodgkin-Huxley model with a slow inward (calcium) current~\cite{Drion_2012}. Both models are planar and lack the second feature, that is, they can only be regarded as a mathematical reduction of a physiological conductance-based model.

This paper aims to contribute to the idea that balancing positive and negative feedback in the slow time scale is a key mechanism to generate rest-spike bistability. This viewpoint is at the core of the planar model in~\cite{Franci_2012} and its importance from a physiological viewpoint is highlighted by~\cite{Franci_2018}. Here we complement those works by studying how this mechanism can be naturally implemented in a physiological context: using two distinct slow currents, one providing negative feedback to restore the membrane potential, the other providing positive feedback to obtain two attractors separated by the stable manifold of a saddle.

The remainder of the paper is organized as follows. \Cref{sec:model} presents the model and recalls the notions of geometric singular perturbation theory needed for its analysis. In \cref{sec:red-dyn} we study numerically the dynamics on the critical manifold, highlighting its persistence properties. \Cref{sec:bistability} builds on this picture to derive conditions for multistability and monostability, we focus on the singular case and mention what hypotheses guarantee persistence. In \cref{sec:variations} we discuss some variations of the same geometric picture, while in \cref{sec:oth-mods} we relate the Hindmarsh-Rose and the transcritical models to the one we are studying. We draw some conclusions in \cref{sec:conclusions}. Two appendices report additional details.

\section{A Model of Rest-Spike Bistability}
\label{sec:model}
We consider a three-dimensional slow-fast conductance-based model defined by
\begin{equation}
  \label{eq:MLslow}
  \begin{aligned}
    \varepsilon \dot v &= i - i_{ion}(v, n, p)\\
    \dot n &= -n + S_n(v) \\
   \tau  \dot p &= - p + S_p(v)
  \end{aligned}
\end{equation}
where $\varepsilon$ is a small parameter. The total ionic current is the sum of  a leak current and three voltage-gated currents:
\begin{equation}
  \begin{aligned}
    i_{ion} &= g_l(v-v_l) + S_m(v)(v - 1) + n(v + 1) + p(v - 1)\\
    &= c(v) + n(v + 1) + p(v - 1)
  \end{aligned}
\end{equation}
The parameters $-1 < v_l < +1$ that appear in the equation can be thought of as reversal potentials. In the absence of an external current $i$, the voltage range $[-1,1]$ is positively invariant. The currents $S_m(v)(v - 1)$ and $p(v - 1)$ are then negative (\emph{inward} currents) whereas the current  $n(v + 1)$ is positive (\emph{outward} current). The variables $n$ and $p$ are gating (positive) variables that model the slow activation of the inward current $p(v - 1)$  and of the outward current  $n(v + 1)$. The inward current $S_m(v)(v - 1)$ has \emph{instantaneous} activation, a standard simplification for currents that activate in the fast time-scale. The functions $S_x(v)$ correspond to activation functions that we assume of the form
\begin{equation}
  \label{eq:activation-function}
  S_x = \frac{g_x}{2} \left(\tanh \left( \frac{v - a_x}{b_x} \right) + 1\right)
\end{equation}
Here the multiplicative factor $g_x$ corresponds to the maximal conductance associated to the current $x$. We find it convenient to write the equations in this form, rather than including maximal conductances in the voltage equation, because this allows us to change maximal conductances of slow currents without modifying the critical manifold of the system. A consequence of this is that the dynamics of $p$ and $n$ lies between zero and the corresponding maximal conductance, i.e. $p \in [0, g_p]$ and $n \in [0, g_n]$.

\begin{figure}[ht]
  \centering
  \extfig{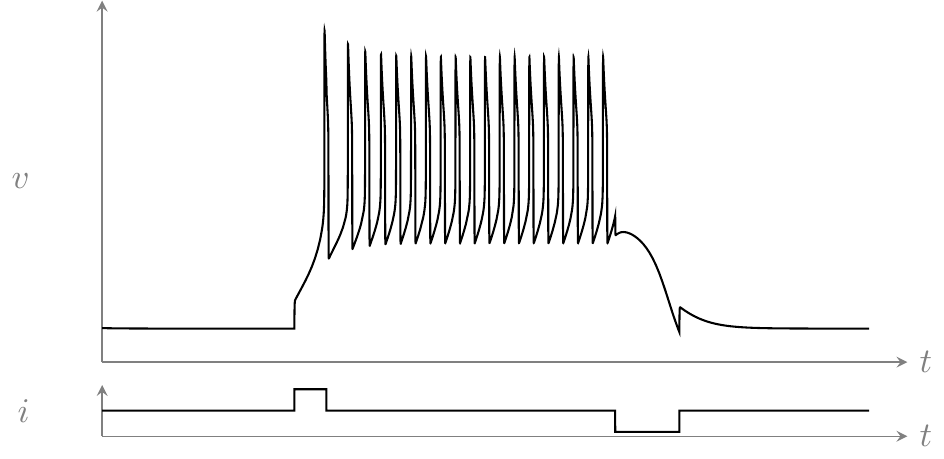}
  \caption{Rest-spike bistability in the model~\cref{eq:MLslow}.}
  \label{fig:rest-spike}
\end{figure}

The key property of the model is the presence of the slow inward current $p(v - 1)$. In the absence of this current, the model is two-dimensional and has a phase portrait similar to the classical FitzHugh-Nagumo model. With this additional slow inward current, both continuous spiking and rest coexist for the same value of applied current, as shown by the simulation in \cref{fig:rest-spike} (see \cref{sec:para} for numerical values of the parameters).

We note that similar phenomena can be obtained with a current of the type $p(v + 1)$ where $p$ inactivates, i.e. decreases as $v$ increases. Physiologically this corresponds to an outward current that inactivates slowly, rather than an inward current that activates slowly. Both types of currents model a source of positive feedback in the slow time-scale~\cite{Franci_2013}. A classical example of  slowly inactivating outward current is the $A$-type potassium current~\cite{Drion_2015}.

We use geometric singular perturbation theory~\cite{Fenichel_1979} to study the slow-fast system~\cref{eq:MLslow} as $\varepsilon$ tends to zero. The singular limit of this model is the differential-algebraic system
\begin{equation}
  \label{eq:slow-dyn}
  \begin{aligned}
    0 &= i - i_{ion}(v, n, p)\\
    \dot n &= -n + S_n(v) \\
    \tau \dot p &= - p + S_p(v)
  \end{aligned}
\end{equation}
which we call slow dynamics or reduced system. After rescaling time, the same limit leads to the layer dynamics
\begin{equation}
  \label{eq:fast-dyn}
  \begin{aligned}
    v' &= i - i_{ion}(v, n, p)\\
    n' &= 0\\
    p' &= 0
  \end{aligned}
\end{equation}
where $'$ refers to differentiation with respect to the fast time $\tau = t / \varepsilon$.

The reduced system~\cref{eq:slow-dyn} is constrained to the critical manifold $\mathcal{C}_0$, defined by
\begin{equation}
  \label{eq:crit_man}
  i_{ion}(v, n, p) = i
\end{equation}
and corresponding to fixed points of the layer dynamics~\cref{eq:fast-dyn}. Normally-hyperbolic compact subsets of $\mathcal{C}_0$ persist as invariant manifolds of~\cref{eq:MLslow} for $\varepsilon$ small enough. This manifolds are not necessarily unique, but we assume one family of perturbation $\mathcal{C}_{\varepsilon}$ has been fixed and call them slow manifolds.

Perturbations of subsets of $\mathcal{C}_0$ maintain their type of stability with corresponding (local) stable and unstable manifolds. These admit invariant foliations, with each point on the critical manifold acting as base for a fiber. Invariance of the foliation can be interpreted as points on each fiber ``shadowing'' the corresponding base point, in forward time for the stable manifold and backward for the unstable. Points on $C_{\varepsilon}$ follow a dynamics that is a regular perturbation of the reduced system~\cref{eq:slow-dyn}; in the following we refer to this perturbation as slow dynamics.

A point $x$ on the critical manifold is normally hyperbolic if it is a hyperbolic fixed point of the layer dynamics~\cref{eq:fast-dyn}. If this is the case, as $\varepsilon\to 0$ the fibers based at $x$ tend to its stable and unstable manifolds in the layer dynamics~\cref{eq:fast-dyn}. For~\cref{eq:MLslow} the layer dynamics is one dimensional, so that hyperbolic fixed points are either attractive or repulsive, with their invariant manifolds corresponding to lines with $n$ and $p$ constant.

\begin{figure}[ht]
  \centering
  \extfig{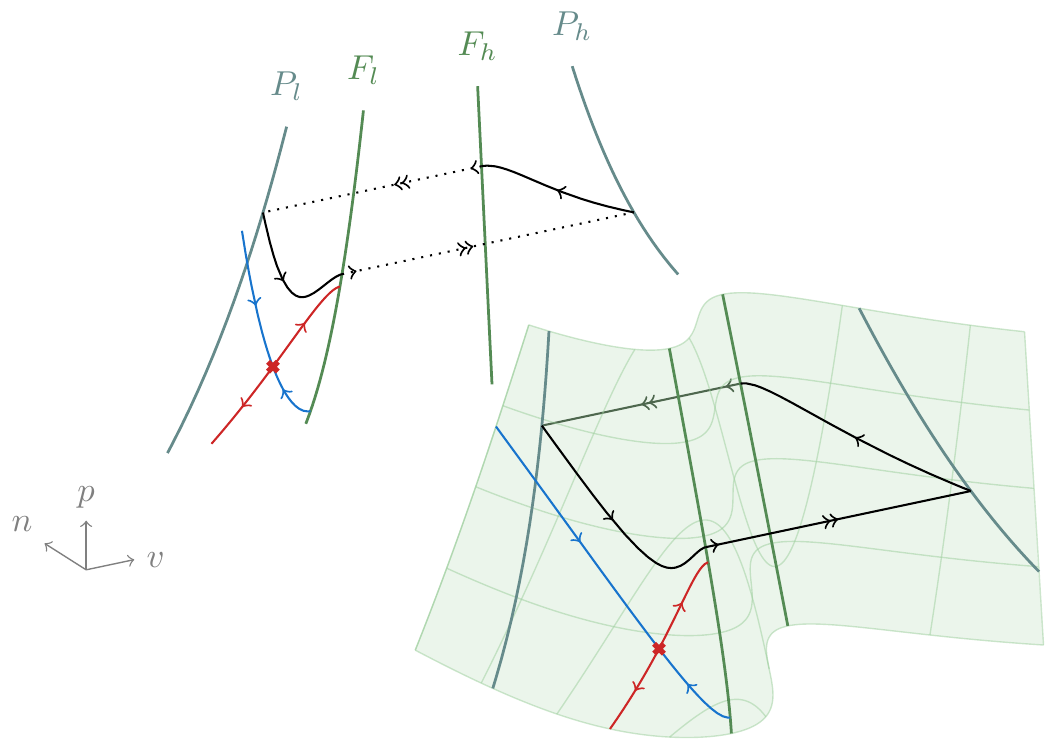}
  \caption{Reduced dynamics~\cref{eq:slow-dyn} on the critical manifold~\cref{eq:crit_man}, and its projection onto the $v$-$p$ plane, together with the lines of folds $F_l$, $F_h$ and their projections $P_l$, $P_h$. A saddle point and its stable (blue) and unstable (red) manifolds in the reduced system are shown on the critical manifold. The black trajectory is a singular relaxation oscillation composed of two slow parts (single arrow) connected by two trajectories along fast fibers (double arrow, dotted in the projection).}
  \label{fig:cman-proj}
\end{figure}

We consider parameter ranges for which the critical manifold can be divided in three normally-hyperbolic branches. These are separated by two lines of folds that we call $F_l$ and $F_h$, and verify
\begin{equation}
  \label{eq:folds}
    \frac{\partial i_{ion}}{\partial v} = \frac{dc}{dv}(v) + n + p = 0
\end{equation}
The two lines of folds are connected by an unstable branch $M$. The other branches, $S_l$ and $S_h$, are both stable. \Cref{fig:cman-proj} shows the typical shape of the critical manifold for fixed $i$.

Around any point away from the lines of folds, the critical manifold admits a parametrization in the slow variables $n$ and $p$. However, this local  parametrization cannot be made global due to the presence of folds. Following~\cite{Szmolyan_2001, Szmolyan_2004, Wechselberger_2013}, we use $v$ and $p$ to obtain a parametrization valid in the interval $v\in(-1,\, 1)$. This is achieved by solving~\cref{eq:crit_man} for $n(v, p, i)$. The corresponding projection is shown in \cref{fig:cman-proj}.

The reduced dynamics in these coordinates is obtained differentiating~\cref{eq:crit_man}:
\begin{equation}
  \label{eq:red-dyn}
  \begin{aligned}
    \frac{\partial i_{ion}}{\partial v} \dot v &= -\frac{\partial i_{ion}}{\partial n} \dot n - \frac{\partial i_{ion}}{\partial p} \dot p\\
    \tau \dot p &= -p + S_p(v)
  \end{aligned}
\end{equation}
The first equation becomes singular on the lines of folds $\frac{\partial i_{ion}}{\partial v} = 0$. Multiplication by $\frac{\partial i_{ion}}{\partial v}$ recovers a regular differential equation
\begin{equation}
  \label{eq:red-dyn-des}
    \begin{aligned}
      \dot v &= -\frac{\partial i_{ion}}{\partial n} \dot n - \frac{\partial i_{ion}}{\partial p} \dot p\\
      \tau \dot p &= \frac{\partial i_{ion}}{\partial v}(-p + S_p(v))
    \end{aligned}
\end{equation}
The two systems~\cref{eq:red-dyn} and~\cref{eq:red-dyn-des} share the same trajectories with different time parametrizations. Moreover, in~\cref{eq:red-dyn-des} time is reversed on the unstable branch $\frac{\partial i_{ion}}{\partial v} < 0$ and new fixed points can appear on the lines of folds. These verify
\begin{equation}
  \label{eq:folded-sing}
    \frac{\partial i_{ion}}{\partial v} = 0 \qquad \frac{\partial i_{ion}}{\partial n} \dot n + \frac{\partial i_{ion}}{\partial p} \dot p = 0
\end{equation}
They are called folded singularities~\cite{Szmolyan_2001}.

Away from the lines of folds the two systems~\cref{eq:red-dyn} and~\cref{eq:red-dyn-des} are largely equivalent, but important differences occur in the neighborhood of $F_l$ and $F_h$. Moreover, near these lines the perturbed dynamics is no longer constrained by normal hyperbolicity, in particular it cannot be obtained as a regular perturbation of the reduced system~\cref{eq:red-dyn}. Different phenomena are possible. The least degenerate situation occurs when the desingularized vector field is never zero along these lines:
\begin{equation}
  \label{eq:jump-point}
  \frac{\partial i_{ion}}{\partial n}\dot n + \frac{\partial i_{ion}}{\partial p} \dot p \not=0
\end{equation}
Under this assumption the desingularized vector field~\cref{eq:red-dyn-des} can point either to the unstable branch or the stable one. Assuming the additional nondegeneracy condition
\begin{equation}
  \label{eq:not-a-cusp}
  \frac{\partial ^2 i_{ion}}{\partial v^2} \not = 0
\end{equation}
the first case corresponds to jump points, at which the reduced system~\cref{eq:red-dyn} admits two solutions backwards in time but none in forward time. For $\varepsilon>0$ a stable branch of $\mathcal{C}_{\varepsilon}$ near these points can be continued using the flow~\cite{Szmolyan_2004}. Doing so shows that trajectories on the slow manifold pass the folds and reach a fiber contained in the stable manifold of the other stable branch of $C_{\varepsilon}$, with the flow contracting the direction transverse to the manifold.

Condition~\cref{eq:jump-point} corresponds to the vector field being transverse to the critical manifold, a condition which  is violated at folded singularities. These are fixed points of the desingularized system~\cref{eq:red-dyn-des}, but not necessarily fixed points of the reduced dynamics~\cref{eq:red-dyn}. As a consequence, they can be reached in finite time. Depending on the type of fixed point they can correspond to the singular limit of canard trajectories, i.e. intersections between stable and unstable branches of the slow manifold~\cite{Szmolyan_2001}. Generically, the desingularized flow changes direction at these points. Hence, a folded singularity delimits the set of jump points on a line of folds~\cite{Wechselberger_2013}.

\section{Reduced Dynamics}
\label{sec:red-dyn}
We will now study the reduced system \cref{eq:red-dyn}, often with the aid of its desingularized version~\cref{eq:red-dyn-des}. Fixed points can be parametrized by $v$ through the steady-state $i$-$v$ curve
\begin{equation}
  \label{eq:i-v-curve}
    i_{s}(v) := i_{ion}(v, S_n(v), S_p(v))
\end{equation}
This is shown in \cref{fig:bifdiag} and is an S-shaped curve, with two folds separating three families of fixed points $\mathcal{X}_l$, $\mathcal{X}_m$ and $\mathcal{X}_h$; $\mathcal{X}_l$ corresponds to low voltages, $\mathcal{X}_m$ to intermediate voltages and $\mathcal{X}_h$ to high voltages. For fixed $i$, we denote points in each family with corresponding lower-case letters $x_l$, $x_m$ and $x_h$. In addition to  these three fixed points, the desingularized dynamics~\cref{eq:red-dyn-des} has a folded singularity $x_f\in F_l$. For parameter values reported in~\cref{sec:para}, and $i$ in the range of interest in this section, this point is a focus and does not lead to canard trajectories~\cite{Szmolyan_2001}; it only delimits jump points on $F_l$.

\begin{figure}[ht]
  \centering
  \extfig{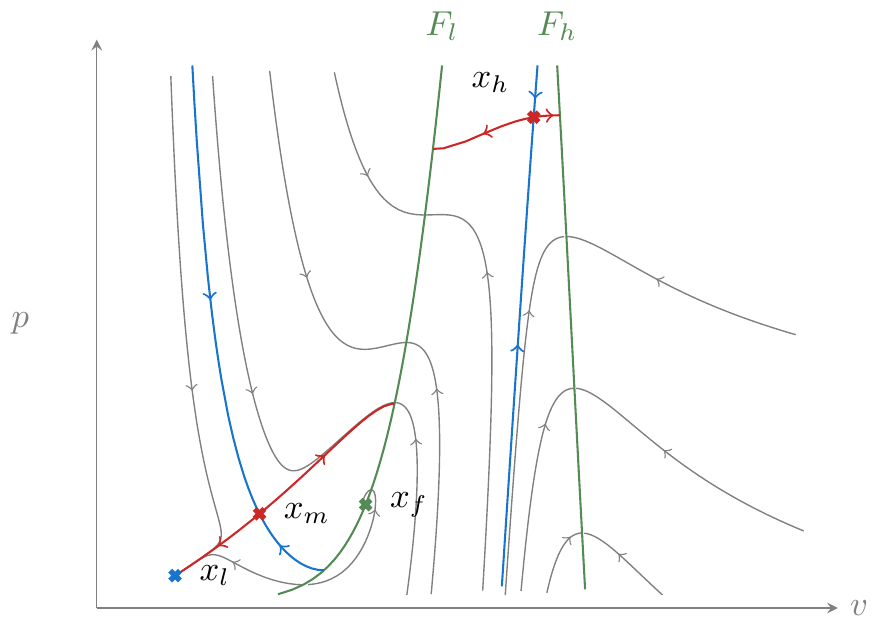}
  \caption{A typical phase portrait of the reduced system~\cref{eq:red-dyn}. Fixed points of the desingularized system~\cref{eq:red-dyn-des} are denoted by crosses, $x_l$ is a stable node, $x_m$ and $x_h$ are saddle points and $x_f$ is a folded focus (unstable). Stable and unstable manifolds of the saddle points are shown in blue and red, respectively. Along the two lines of folds $F_l$ and $F_h$ the system is singular: trajectory at those points are defined only in forward or backward time; the first of these two cases corresponds to jump points. The stable manifold of $x_m$ separates initial conditions in $S_l$ (left of $F_l$) that reach a jump point from those that converge to $x_l$.}
  \label{fig:red-pp}
\end{figure}

\Cref{fig:red-pp} shows the typical phase portrait of the reduced system~\cref{eq:red-dyn}. The fixed point  $x_l$ is a stable node, while $x_m$ and $x_h$ are both saddle points. Their stable and unstable manifolds do not extend beyond $F_l$ and $F_h$ due to loss of existence and uniqueness along these lines. In particular, unstable manifolds terminate at jump points.

For $\varepsilon>0$ hyperbolic fixed points persist in the slow dynamics with their stable and unstable manifolds~\cite{Fenichel_1979,Szmolyan_1991}. In the perturbed system~\cref{eq:MLslow} these fixed points are still hyperbolic. In particular, saddle points remain saddle, with their invariant manifolds being obtained as a combination of trajectories in the slow dynamics and fast fibers. The unstable manifold of $x_m$ is completely contained in the slow manifold. Its stable manifold, instead, is two dimensional; it includes the stable manifold in $\mathcal{C}_{\varepsilon}$ as well as all fast fibers based on that curve. In the singular limit this surface tends to the stable manifold of $x_m$ in the reduced system~\cref{eq:red-dyn} and all nearby segment with constant $n$ and $p$ that intersect it. Similarly, $x_h$ perturbs to a saddle with a one-dimensional stable manifold and a two-dimensional unstable manifold.

Adding a trivial equation for $i$ to~\cref{eq:MLslow}, the same is true for the family of fixed points $x_m(i)$. At least for $i$ in a small interval, this family persists together with its two-dimensional unstable manifold and three-dimensional stable one. Sections of these manifolds for fixed $i$ coincide with the invariant manifolds of the corresponding fixed point.

When $i$ varies on larger domains, the outlined phase portrait can undergo two distinct qualitative changes:
increasing $i$ leads to a fold of the $i$-$v$ curve at which $x_l$ and $x_m$ merge in a saddle-node bifurcation, leaving only one fixed point $x_h\in M$. Likewise, decreasing $i$, $x_m$ and $x_h$ reach a similar fate, leaving $x_l\in S_l$ as the only fixed point.

To obtain this second bifurcation it is necessary that one of the two fixed points crosses the line $F_l$ and changes branch\footnote{This assumes that the bifurcation does not happen exactly on $F_l$.}. In our case $x_m$ crosses $F_l$. This passage corresponds to an exchange of stability with the folded singularity through a folded saddle-node~\cite{Krupa_2010}. Beyond this crossing, the folded singularity is a saddle, while $x_m$ is a node of the reduced system. In a similar fashion increasing the applied current leads to $x_h$ crossing $F_h$, which happens once $x_h$ is the only fixed point left. After this crossing, $x_h$ is a stable fixed point on an attractive branch.

Finally, varying $i$ can lead to changes in the type of folded singularity. As already mentioned $x_m\in F_l$ corresponds to a folded saddle-node, thus varying $i$ and moving $x_m$ between branches leads to different types of folded singularity: it is a saddle when $x_m \in M$ and a node when $x_m \in S_l$. Both situations lead to canard trajectories~\cite{Szmolyan_2001}. Moreover, since $x_f$ is a focus in the phase portrait described above, it has to change to a node before becoming a folded saddle-node.

\section{Rest-spike bistability}
\label{sec:bistability}

Returning to the phase portrait in \cref{fig:red-pp},  we now analyze the global return mechanism that leads to rest-spike bistability.

In the singular limit $\varepsilon = 0$, trajectories on stable branches of the critical manifold $\mathcal{C}_0$ stay on it until they reach a line of fold in correspondence of a jump point. Once one of these points is reached the singular trajectory is continued along a fast fibre with constant $n$ and $p$, reaching the opposite branch as shown in \cref{fig:cman-proj}. The points at which these singular trajectories arrive correspond to the projections of $F_l$ and $F_h$ along fast fibers. We call these projections $P_l \subset S_l$ and $P_h \subset S_h$.

Based on this property, we can analyze the singular system referring only to the $v$-$p$ plane and the reduced dynamics: when a trajectory reaches a jump point it is transported to the corresponding projection keeping $p$ fixed, as shown in \cref{fig:cman-proj} for a limit cycle.



Rest-spike bistability follows from how the stable and unstable manifold of $x_m$ constrain trajectories. The role of the stable manifold is simple, it separates initial conditions on $S_l$ that reach a jump point on $F_l$ from those that remain on the critical manifold and tend to $x_l$. The unstable manifold, instead, determines if the system is multistable. This is the case if the unstable manifold stays away from $x_l$. Otherwise almost all trajectories converge to $x_l$. We treat these two situations separately in the next sections.

\subsection{Bistability}
\label{sec:mul-stab}

\begin{figure}[ht]
  \centering
  \extfig{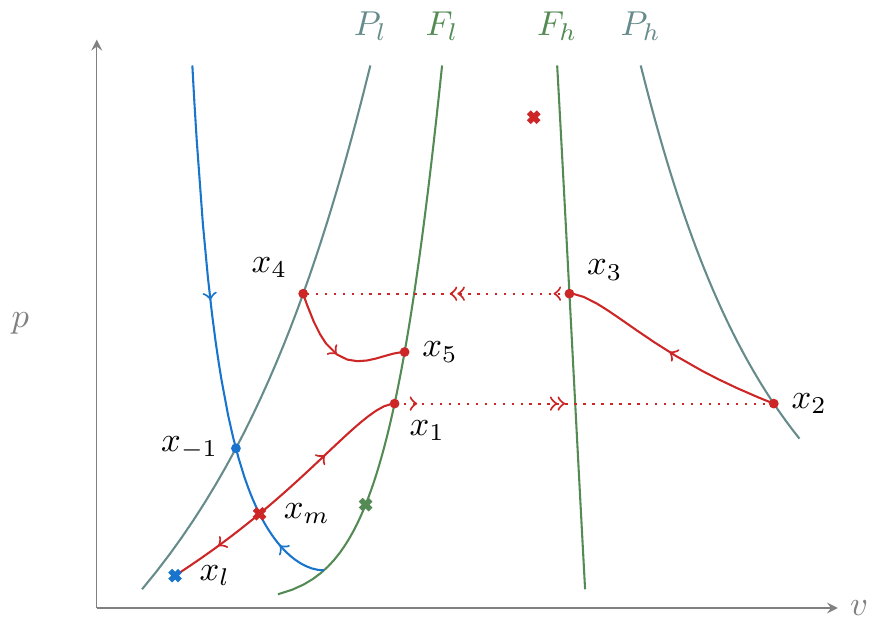}
  \caption{Reduced dynamics~\cref{eq:red-dyn} in the multistable case. The stable manifold of $x_m$ (blue) separates initial conditions that reach a jump point on $F_l$ from those that converge to $x_l$. Jump points are mapped to their projections (e.g. $x_1$ to $x_2$ and $x_3$ to $x_4$). The unstable manifold of $x_m$ (red) delimits an invariant set for the dynamics.}
  \label{fig:pp-limcyc}
\end{figure}

In the following we denote by $x_1$ the intersection of the unstable manifold of $x_m$ with $F_l$, and by $x_{-1}$ the intersection of the stable manifold of $x_m$ with $P_l$. Following the singular flow from $x_1$ leads to $x_2\in P_h$, then to $x_3\in F_h$ and back to $P_l$ at $x_4$ (see \cref{fig:pp-limcyc}). We recall that given the dynamics~\cref{eq:MLslow} we can assume that $p$ lies in the interval $[0, g_p]$, where $g_p$ is the maximal conductance appearing in~\cref{eq:activation-function}.

Assume that the trajectory starting at $x_4$ reaches a jump point on $F_l$ ($x_{5}$), as shown in \cref{fig:pp-limcyc}. Consider the segment $I_l\subset P_l$ between $x_4$ and $p=g_p$. The reduced dynamics maps this segment to $F_l$ in finite time, defining a map $\Pi_l:I_l \to F_l$. Clearly the same map can be defined using the desingularized reduced system~\cref{eq:red-dyn-des}, thus as long as this vector field is transverse to $F_l$ at all points in $\Pi_l(I_l)$ the map is smooth. We note that this is equivalent to $\Pi_l(I_l)$ not containing folded singularities. Similarly, on $S_h$ we define the segment $I_h\subset P_h$ between $x_2$ and $p=g_p$, and a corresponding map $\Pi_h:I_h \to F_h$. We denote the projection along fast fibers by $\Pi_f$ (from $F_l$ to $P_h$ and from $F_h$ to $P_l$). Since the dynamics is bounded by the line $p=g_p$, by construction we have
\begin{equation}
  \label{eq:maps-images}
  \Pi_f \circ \Pi_l (I_l) \subset I_h \qquad \Pi_f \circ \Pi_h(I_h) \subset I_l
\end{equation}
which allows us to define the singular Poincar\'e map
\begin{equation}
  \label{eq:Poincare-map}
  \Pi = \Pi_f \circ \Pi_h \circ \Pi_f \circ \Pi_l: I_l \to I_l
\end{equation}

This construction shows that the stable manifold of $x_m$ divides the state space into two invariant sets. One is the basin of attraction of $x_l$, while the other one has dynamics characterized by the Poincar\'e map~\cref{eq:Poincare-map}. Since this is a smooth map of an interval into itself it admits at least one fixed point, which corresponds to a singular relaxation oscillation. As shown in~\cite{Szmolyan_2004}, if this fixed point is hyperbolic, under the additional hypothesis that the singular trajectory intersects $P_l$ and $P_h$ transversally, it perturbs to a hyperbolic limit cycle for $\varepsilon>0$. In fact the Poincar\'e map~\cref{eq:Poincare-map} is (up to conjugacy) a global version of the one used in that reference.

We remark that this construction only guarantees multistability. Further analysis of the map~\cref{eq:Poincare-map} is required to obtain a more accurate picture. While this is beyond the scope of this work, numerical simulations confirm that this map has a unique attracting fixed point.

\subsection{Monostability}
\label{sec:mono-stab}

\begin{figure}[ht]
  \centering
  \extfig{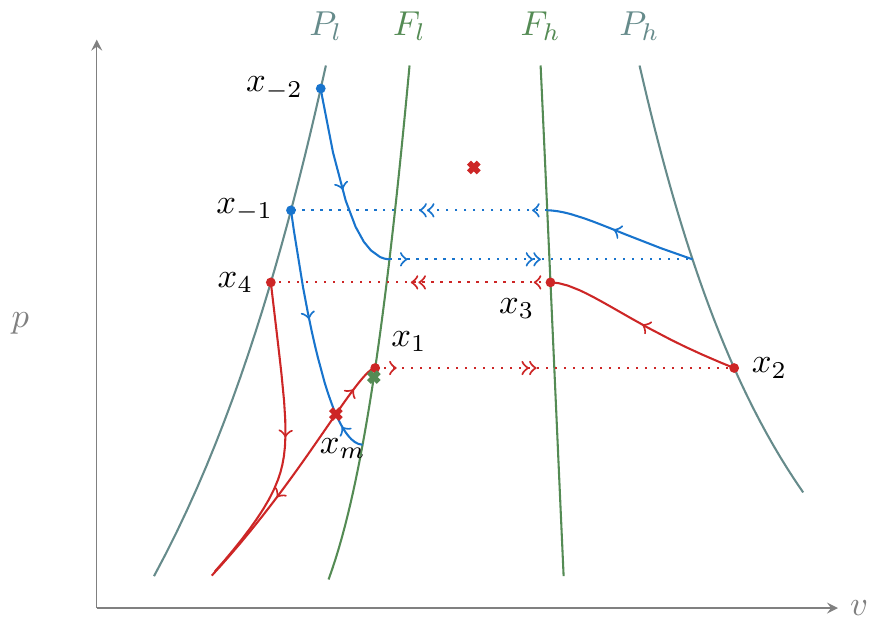}
  \caption{Reduced dynamics in the monostable case. The stable manifold of $x_m$ separates initial condition that arrive to a jump point on $F_l$ from those that converge to $x_l$ (not shown). The unstable manifold of $x_m$ (red) converges to $x_l$ after one jump. Similarly, almost all initial conditions on stable branches converge to it, the only exception being the ones that form the stable manifold of $x_m$.}
  \label{fig:pp-fixpt}
\end{figure}

Constructing the Poincar\'e map~\cref{eq:Poincare-map} requires that $x_4$ falls inside the interval defined by $x_{-1}$ and $p=g_p$ on $P_l$. The situation in which this assumption fails is illustrated  in \cref{fig:pp-fixpt}. In this case most trajectories on $S_l$ and $S_h$ are attracted by the stable fixed point $x_l$, the only exception being the stable manifold of $x_m$.

To see this, we start from $x_{-1}$ and consider its anti-image through $\Pi_f$ on $F_h$. Continuing to follow the singular flow ``backwards'', as shown in \cref{fig:pp-fixpt}, leads back to $P_l$ at a point that we call $x_{-2}$. Any compact segment in $P_l$ that lies between these points is mapped by the singular flow strictly inside the segment delimited by $x_4$ and $x_{-1}$. Since any point strictly inside this second segment converges to $x_l$, the same conclusion extends to all points in the original segment.

The same argument shows that points in the portion of $S_l$ delimited by the trajectories starting at $x_{-1}$ and $x_{-2}$ tend to $x_l$. The only exceptions are these boundary trajectories that reach $x_m$ and belong to its stable manifold. As long as the stable manifold of $x_m$ is unbounded in the $p$ coordinate, the same argument can be iterated on all $S_l$ and adapted to $S_h$, leading to the conclusion that almost all points on $S_l$ and $S_h$ are in the basin of attraction of $x_l$. This situation persists for small enough $\varepsilon>0$, and since most points are attracted to stable branches of the slow manifold we obtain that for almost all initial conditions the perturbed dynamics converges to $x_l$.

\subsection{Homoclinic trajectory and bifurcation diagram}
\label{sec:transition}

Transitions between monostability and bistability in system~\cref{eq:MLslow} are controlled by the applied current $i$. The phase portraits in \cref{fig:pp-limcyc,fig:pp-fixpt} suggest the presence of a homoclinic trajectory, which can be obtained by decreasing  the applied current from the bistable case. In the singular limit this trajectory corresponds to the condition $x_4 = x_{-1}$ and delimits the boundary of bistability. We denote by $i_H$ the value of current at which this happens. While we cannot expect this homoclinic trajectory to persist for $\varepsilon>0$ with $i$ fixed, it is natural to ask whether for $\varepsilon>0$, fixed and small, we can find an $i_{H}(\varepsilon)$, close to $i_H$, at which a homoclinic trajectory exists. There is a natural transversality condition that guarantees this property. The family of fixed points $x_m(i)$ admits a three-dimensional stable manifold and a two-dimensional unstable one. Their intersection is a homoclinic trajectory. In the singular limit, following the unstable manifold of $x_m(i)$ leads back to $S_l$ after two jumps. Extending $\mathcal{C}_0$ to include $i$, $x_m(i)$ is a (normally hyperbolic) invariant set in it, with two-dimensional invariant manifolds. The continuation of the unstable one using the singular flow, after two jumps intersects the stable manifold in the plane $i=i_H$. If this intersection is transverse then it persists for small $\varepsilon$ and $i$ close to $i_H$. We show this in \cref{sec:ex-hom} adapting the arguments used in~\cite{Szmolyan_2004} to prove existence of relaxation oscillations.

\begin{figure}[ht]
  \centering
  \extfig{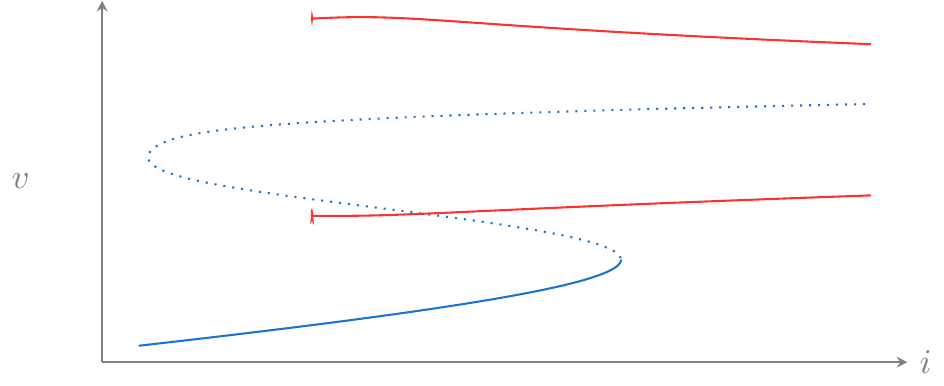}
  \caption{Bifurcation diagram of~\cref{eq:MLslow}. Solid lines denote stable solutions, dotted correspond to unstable ones; blue lines correspond to fixed points, red lines to limit cycle, in the latter case both maximum and minimum are shown.}
  \label{fig:bifdiag}
\end{figure}

To conclude this section, \cref{fig:bifdiag} shows the bifurcation diagram of the whole system~\cref{eq:MLslow} computed with AUTO-07p~\cite{doedel2007auto} for parameter values reported in~\cref{sec:para}. The numerics confirms the presence of a family of limit cycle (red curves) and its coexistence with a family of fixed points (blue curve). The family of periodic solutions terminates in a homoclinic trajectory for low values of $i$ (the numerical continuation was stopped at period $T=10^4$).

\section{A common geometric picture}
\label{sec:variations}
The bifurcation diagram illustrated in the previous section is understandably only one among many possible scenarios compatible with the three-dimensional geometry of \cref{fig:cman-proj}.
While a detailed study of all possible cases is beyond the scope of this work, we wish to  highlight how different types of bistability could share the same geometric structure. To do this we use ideas and techniques from~\cite{ribar2018neuromorphic}. As in~\cref{sec:red-dyn} we identify fixed points with the $i$-$v$ curve
\begin{equation}
  \label{eq:iss}
    i_{s}(v) = i_{ion}(v, S_n(v), S_p(v))
\end{equation}
and divide them in three families $\mathcal{X}_l$, $\mathcal{X}_m$ and $\mathcal{X}_h$, separated by two folds. As noted in~\cref{sec:red-dyn} there is a value of current $i_c$, between the two folds, at which $x_m$ crosses $F_l$ to enter the unstable branch $M$. The scenario studied in~\cref{sec:bistability} assumes  $i_c < i_H$ since the homoclinic bifurcation occurs when $x_m \in S_l$.

\begin{figure}[ht]
  \centering
  \extfig{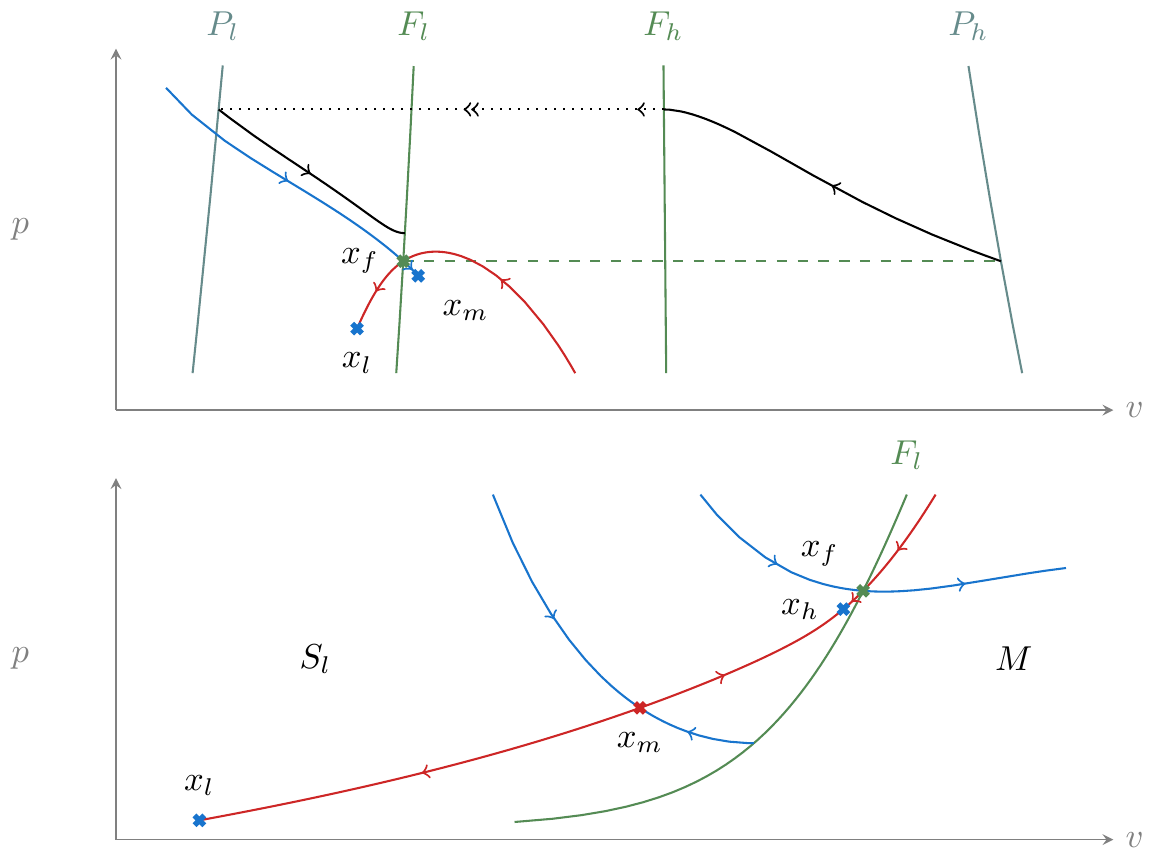}
  \caption{Alternative scenarios that lead to bistability. Top: geometric construction when a folded saddle ($x_f$) takes the place of $x_m$. Bottom: bistability between two fixed points ($x_l$ and $x_h$).}
  \label{fig:pp-others}
\end{figure}

As a first variation we consider what happens when the bistable range extends to current values for which  $x_m\in M$. The bifurcation $x_m\in F_l$ corresponds to a folded saddle-node. Beyond this  bifurcation $x_m\in M$ is a node of the reduced dynamics while $x_f$ is a saddle. In this case the analysis is easily adapted from \cref{sec:bistability}. One must simply  substitute the stable manifold of $x_m$ with the one of $x_f$, and use $\Pi_f(x_f)$ in place of $x_2 = \Pi_f(x_1)$. \Cref{fig:pp-others} shows the corresponding geometric construction. A classical example where this scenario occurs is the Hodgkin-Huxley model with the reversal potential of potassium increased. This situation of  bistability has been studied in the early work~\cite{rinzel1985excitation}. Its planar reduction leads to the transcritical model~\cite{Franci_2013}. Also in this case the boundary of bistability is a singular homoclinic trajectory. This trajectory, however, has to go through the folded singularity $x_f$ to reach $x_m$ on the unstable branch $M$.

Both cases discussed so far assume that $x_m$ and $x_h$ collide in a fold on $M$. Yet another scenario corresponds to this fold occurring on $S_l$, after $x_h$ crosses $F_l$. Also this crossing leads to a folded saddle-node, after which $x_h\in S_l$ can perturb to a stable fixed point. Local analysis around folded saddle-node shows the possibility of Hopf bifurcations~\cite{Krupa_2010}, which are indeed  found numerically. After this the system presents two stable fixed points. The relevant part of the reduced dynamics in this case is shown in \cref{fig:pp-others}: the stable manifold of $x_m$ acts as separatrix between the basins of attraction of the two stable fixed points, while the one of $x_f$ (a folded saddle) separates initial conditions that reach a jump point on $F_l$ from those that remain on the critical manifold.

The examples above suggest that many possible variants for transitions between monostability and bistability are possible. We also note that many of the geometric constructions used in~\cite{Franci_2012,Franci_2018,ribar2018neuromorphic} have an analog in our setting, allowing, for example, non-plateau oscillations, contrary to the case showed in \cref{fig:rest-spike}. This flexibility is interesting in the perspective of connecting the present approach to the classification of bursting types according to the transitions that occur from rest to spike and vice-versa  (see e.g. \cite{IZHIKEVICH_2000}).

\section{Connections with phase portrait analysis}
\label{sec:oth-mods}

We close this paper by clarifying the connection between the proposed three-dimensional model and two published slow-fast phase portraits of rest-spike bistability.

\begin{figure}[!ht]
  \centering
  \extfig{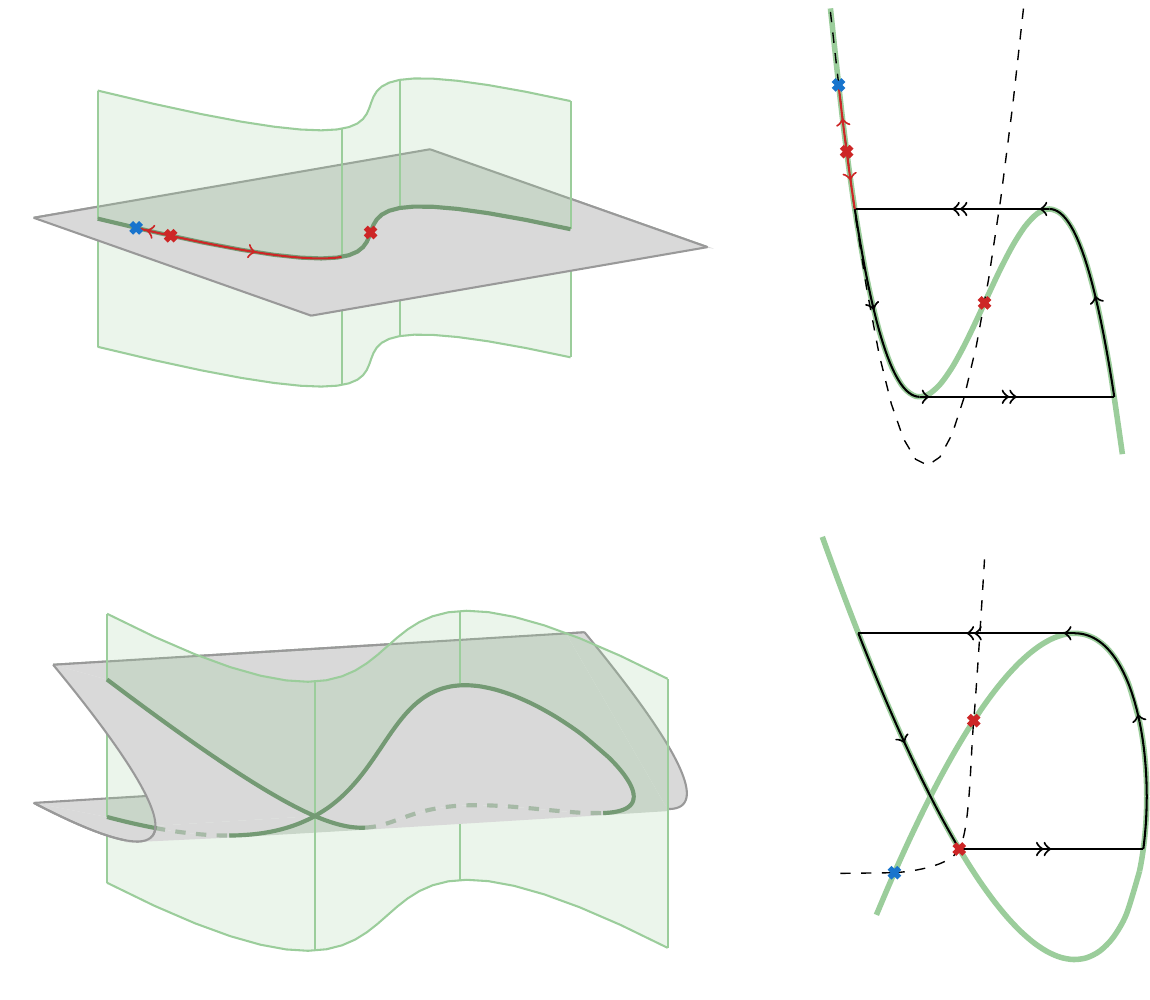}
  \caption{Bistable slow-fast phase portraits as reduction of a larger dimensional model. Left: critical manifolds obtained as the intersection of a higher-dimensional one (green) with a surface (gray). Right: corresponding phase plane with the critical manifold obtained (green) and a possible nullcline for the slow variable (dashed) that completes the dynamics. Top: Hindmarsh-Rose model can be obtained constraining the dynamics to a plane, the critical manifold in the phase plane is the classical N-shaped one, but presents nontrivial dynamics leading to rest-spike bistability. Bottom: the transicritical model obtained constraining the dynamics to a surface. The transcritical bifurcation is obtained when this surface is tangent to a line of folds at a point. This bifurcation is responsible for a singular homoclinic trajectory in the planar reduction.}
  \label{fig:planar}
\end{figure}

The first phase portrait goes back to the seminal work of Hindmarsh and Rose~\cite{Hindmarsh_1982,Hindmarsh_1984}. In one of the earliest attempts to model slow spiking and bursting, Hindmarsh and Rose proposed to modify FitzHugh-Nagumo model with a recovery variable that has a nonmonotonic activation function. Geometrically, this situation corresponds to a degenerate case of the planar pictures described in \cref{sec:red-dyn} and \cref{sec:bistability}, in which all essential elements are contained on a line. As a result, the main elements of the three dimensional dynamics can be captured by constraining it to a plane, resulting in a simplified two-dimensional model of rest-spike bistability. This is characterized by the classical N-shaped critical manifold, as shown in \cref{fig:planar}. The price paid for this simplification is that the flexibility of the two-dimensional slow dynamics described in \cref{sec:variations} is lost. For instance, bistability is only possible if $x_l$ lies out of the stripe delimited by $P_l$ and $F_l$, ruling out patterns in which the voltage of the resting state is between maximum and minimum of the spike.  We note that the nonmonotonicity of the activation function in Hindmarsh-Rose model has the natural interpretation of summarising in one variable the distinct roles of an inward and an outward slow current.

The second rest-spike bistable phase portrait  is the transcritical model of~\cite{Franci_2012}. This model was obtained as a two-dimensional reduction of a conductance-based model that adds a slow calcium current to the Hodgkin-Huxley model~\cite{Drion_2012}. The analysis of~\cite{Franci_2012} rests on the presence of a transcritical bifurcation of the critical manifold. This bifurcation also directly relates to the mixed role of the slow variable as a source of both positive and negative feedback in the slow time-scale. A main motivation of the present paper was to understand the geometric picture generated by this motif in conductance-base models, where these two roles are often played by distinct variables.

To connect the transcritical bifurcation of the planar model~\cite{Franci_2012} to the three-dimensional geometry of the present paper we consider how this planar reduction can be obtained. Referring to our model~\cref{eq:MLslow} for simplicity, a planar reduction is typically obtained imposing an algebraic constraint between $n$ and $p$, which can be interpreted as a path $n(s)$, $p(s)$~\cite{rinzel1985excitation}. After obtaining a dynamic equation for $s$ from a combination of $\dot n$ and $\dot p$, the system becomes
\begin{equation}
  \label{eq:planar-reduction}
  \begin{aligned}
    \varepsilon \dot v &= i - i_{ion}(v, n(s), p(s))\\
    \dot s &= g(v, s)
  \end{aligned}
\end{equation}
which is a slow-fast planar model. Its critical manifold is given by
\begin{equation}
  \label{eq:crit-man-red-mod}
    i = i_{ion}(v, n(s), p(s))
\end{equation}
It corresponds to the intersection of the critical manifold of the larger system with the surface
\begin{equation}
  \label{eq:red-surf}
    n = n(s) \qquad p = p(s)
\end{equation}
A transcritical bifurcation is obtained when
\begin{equation}
  \label{eq:trans-conds}
  \begin{aligned}
    i_{ion} &= i\\
    \frac{\partial i_{ion}}{\partial v} &= 0 \\
    \frac{\partial}{\partial s}(i_{ion}(v, n(s), p(s))) &= \frac{\partial i_{ion}}{\partial n} \frac{d n}{d s} + \frac{\partial i_{ion}}{\partial p} \frac{d p}{d s}=0
  \end{aligned}
\end{equation}
Geometrically this corresponds to a point at which the surface~\cref{eq:red-surf} is tangent to the line of folds of the critical manifold, as shown in \cref{fig:planar}. Similar geometric constructions lead to the presence of a transcritical bifurcation when reducing the Hodgkin-Huxley model with increased potassium reversal potential, as well as when reducing the same model augmented with a calcium current, as done in~\cite{Franci_2012}.

An equivalent interpretation of how the transcritical bifurcation arises is that the path $(n(s), \, p(s))$ defining the surface~\cref{eq:red-surf} is tangent to the line of folds
\begin{equation}
  \label{eq:folds2}
    i = i_{ion}(v, n, p) \qquad \frac{\partial i_{ion}}{\partial v}(v, n, p) = 0
\end{equation}
projected onto the $n$-$p$ plane. This is the simplest example of how singularities in the sense of~\cite{Golubitsky_1985} can be generated from elementary catastrophes, the core idea in the path formulation of~\cite[Ch.3~\S 12]{Golubitsky_1985}. This is particularly interesting in view of~\cite{Franci_2014}, where singularity theory is used to obtain a global description of the critical manifolds of slow-fast planar systems relevant to neuronal dynamics. Two singularities play a prominent role: hysteresis, in connection with spiking, and winged cusp, for rest-spike bistability. Both these singularities can be realized as paths in the unfolding of the cusp catastrophe~\cite{Golubitsky_1985}. Interestingly, this bifurcation is often found in the fast subsystem of neuronal models (an early example being~\cite{ZEEMAN_1973}), and it is typically related to the appearance and disappearance of bistability. For example, decreasing the sodium conductance in the Hodgkin-Huxley model leads to the appearance of this bifurcation, and the same is achieved by reducing $g_m$ in~\cref{eq:MLslow}. The presence of this type of bifurcation in these models suggests that those singularities can arise from model reduction similarly to what happens in the transcritical case.

\section{Conclusions}
\label{sec:conclusions}

We studied a simplified slow-fast model of neuronal activity that exhibits rest-spike bistability. The simplest physiological models of excitability include a fast-activating inward current and a slowly-activating outward current. Our model adds a slowly-activating inward current to this basic motif. We think of this model as a core structure for the generation of multistability in more general and realistic conductance-based models. We speculate that similar results are possible using a slowly inactivating outward current, which would have the same functional role of a slow positive feedback.

Through geometric singular perturbation theory we could analyze the geometry of this three-dimensional model. This geometry is rather simple, with the slow dynamics taking place on a classical N-shaped critical manifold. The saddle point on the critical manifold is a key feature of the proposed model. Its stable manifold acts as separatrix, while its unstable manifold determines whether multiple attractors are present. Moreover, a same geometric picture captures different types of bistability, suggesting a common framework to study different phenomena important to neuronal dynamics.

This is by no means the first study of a slow-fast systems with one fast and two slow variables, nor the first single-cell model of bistability. The value of this model is in that it explains how bistability can arise in a physiologically relevant context using a mechanism that is generic but not widely acknowledged. Our hope is that it contributes to the view that a combination of positive and negative feedback in the slow time-scale is a core element in the generation of neuronal patterns.

\appendix
\section{Existence of a homoclinic trajectory}
\label{sec:ex-hom}
In this section we show that under the assumption of transversality, the intersection of stable and unstable manifolds that leads to a singular homoclinic trajectory persists for $\varepsilon>0$. We do this using the setting of~\cite{Szmolyan_2004} and in particular their results on maps defined by the flow of~\cref{eq:MLslow}. We recast these results in the notation of~\cref{sec:bistability} and refer the reader to the original work for details.

As in \cref{sec:bistability}, $i_H$ is the value of $i$ at which a singular homoclinic trajectory exists. We consider the reduced dynamics for this value of $i$, and fix a point $x_u$ on the unstable manifold of $x_m$ between $x_m$ and $F_l$. Similarly, we fix a point $x_s$ on the stable manifold between $P_l$ and $x_m$.

After a local change of coordinates we can find two neighborhoods of these points, $N_s$ and $N_u$, such that the critical manifold $\mathcal{C}_0$ corresponds to the plane $v=0$. The intersections of these neighborhoods with the planes $n=n_s$ and $n=n_u$ determine two surfaces $\Sigma_s$ and $\Sigma_u$. Rotating $n$ and $p$ if necessary, we can assume that $\Sigma_u \cap \mathcal{C}_0$ intersects the unstable manifold of $x_m$ transversally and only at $x_u$, and similarly for $\Sigma_s \cap \mathcal{C}_0$. For fixed $\delta$ we let $N_{\delta} = (i_H - \delta, \, i_H + \delta)$ %
and consider $\Sigma_s \times N_{\delta}$ and $\Sigma_u \times N_{\delta}$. If $\delta$ is small enough, stable and unstable manifolds of $x_m(i)$ intersect transversally these extended neighborhoods (in the critical manifold extended to include $i$). In the following we assume that $N_s$, $N_u$, $N_{\delta}$ are shrunk whenever necessary.

In~\cref{sec:red-dyn}, we have characterized the stable manifold of $x_m$ for small $\varepsilon >0$, this is composed of a line on $\mathcal{C}_{\varepsilon}$ and the fibers based on it. In the limit $\varepsilon\to 0$, the singular stable manifold intersect $\Sigma_s$ transversally along one of these fibers. Thus if $\varepsilon$ and $\delta$ are small enough the same will be true for the stable manifold of $x_m(i)$ for fixed $i$ and $\varepsilon$. Moreover, since at $\varepsilon=0$, $i=i_H$ this intersection is a line of constant $p$, we can find a parametrization of it that has the form $p = p_s(v, i, \varepsilon)$. Similarly, the intersection of $\Sigma_u$ with the unstable manifold of $x_m(i, \varepsilon)$ defines two functions $v_u(i, \varepsilon)$ and $p_u(i, \varepsilon)$.

Notice that in this section we use $v$ and $p$ to parametrize the two slices $\Sigma_s$ and $\Sigma_u$, so that $v$ preserves its nature of fast variable. This differs from the use of $v$ and $p$ to parametrize the critical manifold as done in~\cref{sec:red-dyn} and~\cref{sec:bistability}.

We can now use the same construction of~\cite{Szmolyan_2004} to obtain a map $\Pi:\Sigma_u \to \Sigma_s$ corresponding to the action of the flow. This has the form
\begin{equation}
  \label{eq:Pi-map}
  \Pi \begin{pmatrix}
    v \\ p
\end{pmatrix} = \begin{pmatrix}
  R(v, p, i, \varepsilon)\\
  G(v, p, i, \varepsilon)
\end{pmatrix}
\end{equation}
$R$ is exponentially small in $\varepsilon$ ($\lvert R \rvert + \lVert \nabla R \rVert < \exp(-c/\varepsilon)$) and in particular verifies
\begin{equation}
  \label{eq:R-at-zero}
  R(v, p, i, 0) = 0
\end{equation}
$G$ has the form
\begin{equation}
  \label{eq:G-form}
  G = G_0(p) + \mathcal{O}(\varepsilon \ln(\varepsilon))
\end{equation}
where $G_0:\Sigma_u \cap \mathcal{C}_0 \to \Sigma_s \cap \mathcal{C}_0$ is the map defined by the singular flow. Smooth dependence on $i$ follows from standard results.

The only difference between this map and the Poincar\'e map defined in~\cite{Szmolyan_2004} is that we consider two different sections $\Sigma_s$ and $\Sigma_u$ rather than one.

Applying this map to $(v_u, \,p_u)$ we obtain the intersection of the unstable manifold of $x_m$ with $\Sigma_s$
\begin{equation}
  \label{eq:Pi-xu}
    \Pi \begin{pmatrix}
      v_u \\ p_u
  \end{pmatrix} = \begin{pmatrix}
    R(v_u, p_u, i, \varepsilon)\\
    G(v_u, p_u, i, \varepsilon)
  \end{pmatrix}
\end{equation}
In this setting an intersection of stable and unstable manifolds corresponds to a solutions of
\begin{equation}
  \label{eq:intersect}
  G(v_u, p_u, i, \varepsilon) = p_s(R(v_u, p_u, i, \varepsilon), i, \varepsilon)
\end{equation}
where $v_u=v_u(i, \varepsilon)$ and $p_u=p_u(i, \varepsilon)$. Thus, we can define
\begin{equation}
  \label{eq:Pu-Ps}
    P_u(i, \varepsilon) = G(v_u, p_u, i, \varepsilon) \qquad P_s(i, \varepsilon) = p_s(R(v_u, p_u, i, \varepsilon), i, \varepsilon)
\end{equation}
and a homoclinic trajectory corresponds to $P_u - P_s=0$.

At $\varepsilon=0$
\begin{equation}
  \label{eq:Pu-Ps-eps0}
  \begin{aligned}
    P_u(i, 0) &= G_0(p_u(i, 0))\\
    P_s(i, 0) &= p_s(R(0, p_u(i, 0), i, 0), i, 0) = p_s(0, i, 0)
  \end{aligned}
\end{equation}
and the existence of the singular homoclinic trajectory at $i = i_H$ means that
\begin{equation}
  \label{eq:homo-cond}
  P_u(i_H, 0) = G_0(p_u(i_H, 0)) = p_s(0, i_H, 0) = P_s(i_H, 0)
\end{equation}
Assuming that
\begin{equation}
  \label{eq:transv-cond}
    \frac{\partial P_u}{\partial i}(i_H, 0) - \frac{\partial P_s}{\partial i}(i_H, 0) \not = 0
\end{equation}
an application of the implicit function theorem\footnote{We need a relaxed version of the implicit function theorem since the dependence on $\varepsilon$ is not smooth but only continuous, this is proven in~\cite{Loomis_2013}.} guarantees the existence of a continuous functions $i_H(\varepsilon)$ such that $P_u(i_H(\varepsilon), \varepsilon) = P_s(i_H(\varepsilon), \varepsilon)$.

At $\varepsilon=0$, $P_u(i, 0)$ is the intersection of the singular unstable manifold (after two jumps) with $\Sigma_s \cap \mathcal{C}_0$, while $P_s(i, 0)$ corresponds to the intersection of the stable manifold of the reduced flow and $\Sigma_s \cap \mathcal{C}_0$. Condition~\cref{eq:transv-cond} corresponds to transversality of the intersection between the manifolds $p=P_s(i,0)$ and $p=P_u(i,0)$ in the extended neighborhood $\Sigma_s \times N_{\delta}$. Since the invariant manifolds of $x_m(i)$ can be obtained applying the singular flow to these two sections, we see that condition~\cref{eq:transv-cond} is equivalent to transversality of the intersection between the invariant manifolds of $x_m(i)$ on the critical manifold (where the unstable manifold has been continued past two jumps using the singular flow).

\section{Parameters}
\label{sec:para}
The analysis in \cref{sec:red-dyn,sec:bistability} uses the following numerical values for the parameters of~\cref{eq:MLslow}
\begin{equation}
  \label{eq:para}
  \begin{gathered}
    \varepsilon = 0.05, \quad v_l = -.8, \quad gl = 2 \\
    g_m = 4.4, \quad a_m = -.19, \quad b_m = .18 \\
    g_n = 8.0, \quad a_n = -.16, \quad b_n = .29 \\
    g_p = 2.0, \quad a_p = -.5, \quad b_p = .3 \\
    \tau = 1.5
  \end{gathered}
\end{equation}
Parameters related to $m$, $n$ and leak are approximated from those in~\cite{Ermentrout_2010} for Morris-Lecar model in the case Hopf, after normalizing $v$ between -1 and 1. Other parameters have been chosen to obtain the desired bistability.

\bibliographystyle{siamplain}
\bibliography{rest-spike}
\end{document}